\documentclass{article}
\usepackage{amsmath,amssymb,latexsym, epsfig, multicol,tikz}
\newtheorem{lemma}{Lemma}
\newtheorem{theorem}{Theorem}
\newtheorem{corollary}{Corollary}

\newcommand{\Cayley}{{\rm Cayley}}

\title{Ricci-flat graphs with girth at least five}
\author{Yong Lin \\ Renmin University of China
\thanks{Renmin University of China, Beijing 100872, China,({\tt linyong01@ruc.edu.cn}).
Supported by National Natural Science
Foundation of China (Grant Nos 11271011). } \and Linyuan Lu \\
University of South Carolina
\thanks{University of South Carolina, Columbia, SC 29208,
({\tt lu@math.sc.edu}).
Supported in part by NSF grant
DMS 1000475. }
\and S.-T. Yau \\Harvard University
\thanks{Harvard University, Cambridge, MA 02138,({\tt yau@math.harvard.edu}).}
}

\begin{document}
\maketitle

\begin{abstract}
  A graph is called {\em Ricci-flat} if its Ricci-curvatures vanish on all
  edges. Here we use the definition of Ricci-cruvature on graphs given
  in [Lin-Lu-Yau, Tohoku Math., 2011], which is a variation of
  [Ollivier, J. Funct. Math., 2009].  In this paper, we classified all
  Ricci-flat connected graphs with girth at least five: they are the infinite
  path, cycle $C_n$ ($n\geq 6$), the dodecahedral graph, the Petersen graph,
   and the half-dodecahedral graph. We also construct many Ricci-flat
   graphs with girth $3$ or $4$ by using the root systems of simple
   Lie algebras.
\end{abstract}

\section{Introduction}

The Ricci curvature plays a very important role in geometric analysis
on Riemannian manifolds. Ricci-flat manifolds are Riemannian maniflods
with Ricci curvature vanishes.  In physics, they represent vacuum
solutions to the analogues of Einstein's equations for Riemannian
manifolds of any dimension, with vanishing cosmological constant. The
important class of Ricci-flat manifolds is Calabi-Yau manifolds.  This
follows from Yau's proof of the Calabi conjecture, which implies that
a compact $K\ddot{a}hler$ manifold with a vanishing first real Chern
class has a $K\ddot{a}hler$ metric in the same class with vanishing
Ricci curvature. They are many works to find the Calabi-Yau manifolds.
Yau conjectured that there are finitely many topological types of
compact Calabi-Yau manifolds in each dimension. This conjecture is
still open.  In this paper, we will study this question on
graphs. First we will give a short history of the definition of Ricci
curvature in discrete setting.

The definition of the Ricci curvature on metric spaces was first from
the well-known Bakry and Emery notation. Bakry and Emery\cite{be1985}
found a way to define the ``lower Ricci curvature bound'' through the
heat semigroup $(P_t)_{ t\ge 0}$ on a metric measure space $M$. There
are some recent works on giving a good notion for a metric measure
space to have a ``lower Ricci curvature bound'', see
\cite{Sturm-2006}, \cite{Lott-2006} and \cite{Ohta-2006}. Those
notations of Ricci curvature work on so-called length spaces. In 2009,
Ollivier \cite{ollivier} gave a notion of coarse Ricci curvature of
Markov chains valid on arbitrary metric spaces, such as graphs.

Graphs and manifolds are quite different in their nature. But they
do share some similar properties through Laplace operators, heat
kernels, and random walks, etc. Many pioneering works were done by
Chung, Yau, and their coauthors \cite{cy95, cy95', cgy96, cy96,
cy97, cgy97, cy99, cy99', cy00, cy00', cgy00}.

A graph $G=(V,E)$ is a pair of the vertex-set $V$ and the edge-set
$E$. Each edge is an unordered pair of two vertices. Unless
otherwise specified, we always assume a graph $G$ is simple (no
loops and no multi-edges) and connected. It may have infinite but
countable number of vertices.  For each vertex $v$, the degree $d_v$
is always bounded. Starting from a vertex $v_1$ we select a vertex
$v_2$ in the neighborhood of $v_1$ at random and move to $v_2$
then we select a vertex $v_3$ in the neighborhood of $v_2$
at random and move to $v_3$, etc. The random sequence of
vertices selected this way is a random walk on the graph. Ollivier
\cite{ollivier}'s definition of the coarse Ricci curvature of Markov
chains on metric space can be naturally defined over such graphs.

The first definition of Ricci curvature on graphs was introduced by
Fan Chung and Yau in 1996 \cite{cy96}. In the course of obtaining a
good log-Sobolev inequality, they found the following definition of
Ricci curvature to be useful:

We say that a regular graph $G$ has a local $k$-frame at a vertex
$x$ if there exist injective mappings $\eta_1, \ldots, \eta_k$ from
a neighborhood of $x$ into $V$ so that

(1)  $x$ is adjacent to $\eta_i x$ for $1\leq i \leq k\,$;

(2) $\eta_i \,x  \neq \eta_j \,x~~$ if $i \neq j\,$.

The graph $G$ is said to be Ricci-flat at $x$ if there is a local
$k$-frame in a neighborhood of $x$ so that for all $i\,$,
$$\bigcup_j \left(\eta_i \eta_j\right) x = \bigcup_j \left(\eta_j \eta_i\right) x~.$$

For a more general definition of Ricci curvature, in \cite{ly09},
Lin and Yau  give a generalization of lower Ricci curvature bound
in the framework of graphs in term the notation of Bakry and Emery.

In our previous paper \cite{LLY}, the Ricci curvature on graphs is
defined based on Ollivier's definition of Ricci curvature for Markov
chains on graphs.  It is natural to define a {\em Ricci-flat graph} to be
a graph where Ricci-curvature vanishes on every edge. This definition 
does not require a graph to be regular; which is an advantage over
the Chung-Yau's definition. The Ricci flat
graphs defined by Chung and Yau are not necessarily Ricci-flat in the
sense of our definition. However, the Ricci curvatures of those graphs 
are always non-negative. In the last section, we constructed many
``Ricci-flat'' graphs under both definitions.

A well-known Bonnet-Myers theorem on Riemannian geometry said that if a
complete Riemannian manifolds with Ricci curvature bounded below by a
positive constant, then it is compact and has a finite fundamental
group.  In the paper of \cite{LLY}, we prove the first part result
of Bonnet-Myers theorem on graphs with Ricci curvature bounded below
by a positive constant. 
In the paper of \cite{forman}, Forman introduced the Ricci curvature
on cell complexes and also obtained the Myers theorem on the fundamental
group of the complexes.

In  this paper, we classified 
Ricci flat graphs with large girth (using our definition).
\begin{theorem}\label{t1}
Suppose that $G$ is a Ricci flat graph with girth $g(G)\geq 5$. Then
$G$ is one of the following graphs,
\begin{enumerate}
\item the infinite path,
\item  cycle $C_n$ with $n\geq 6$,
\item the  dodecahedral graph,
\item the Petersen graph,
\item the half-dodecahedral graph.
\end{enumerate}
\end{theorem}

\begin{figure}[htbp]
\centerline{ {\psfig{figure=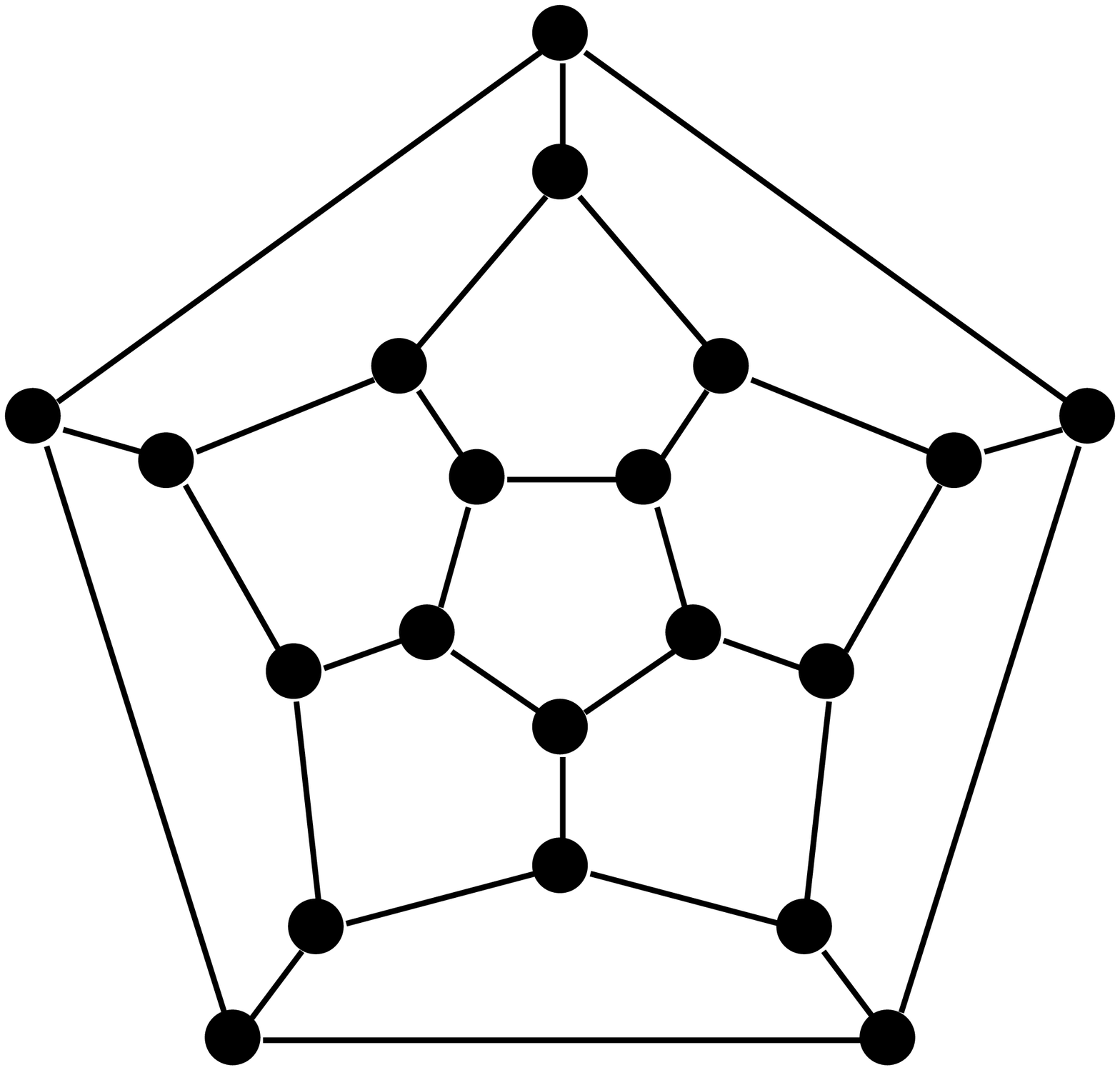, width=0.3\textwidth}}
\hfil \psfig{figure=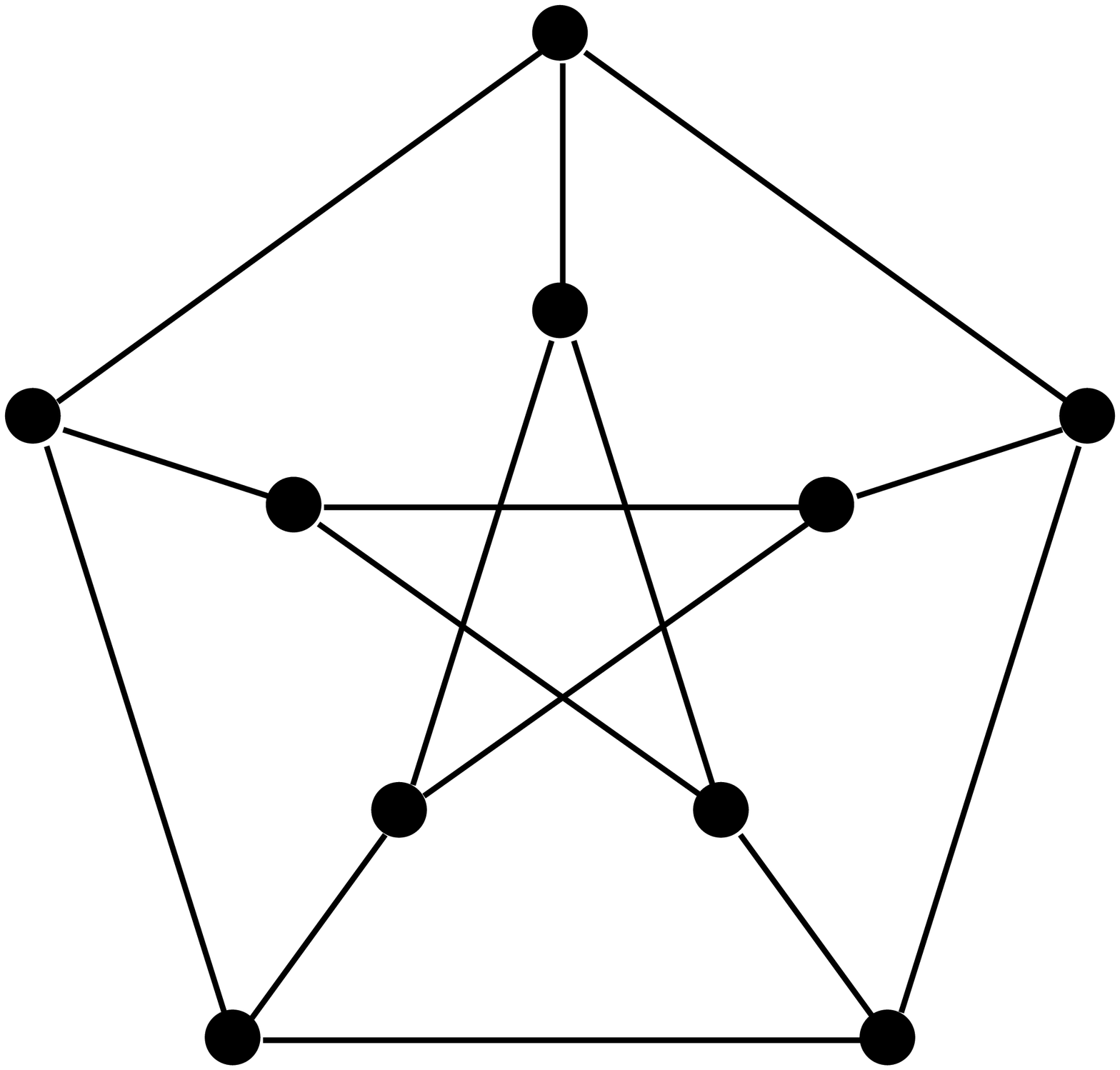, width=0.3\textwidth}
\hfil {\psfig{figure=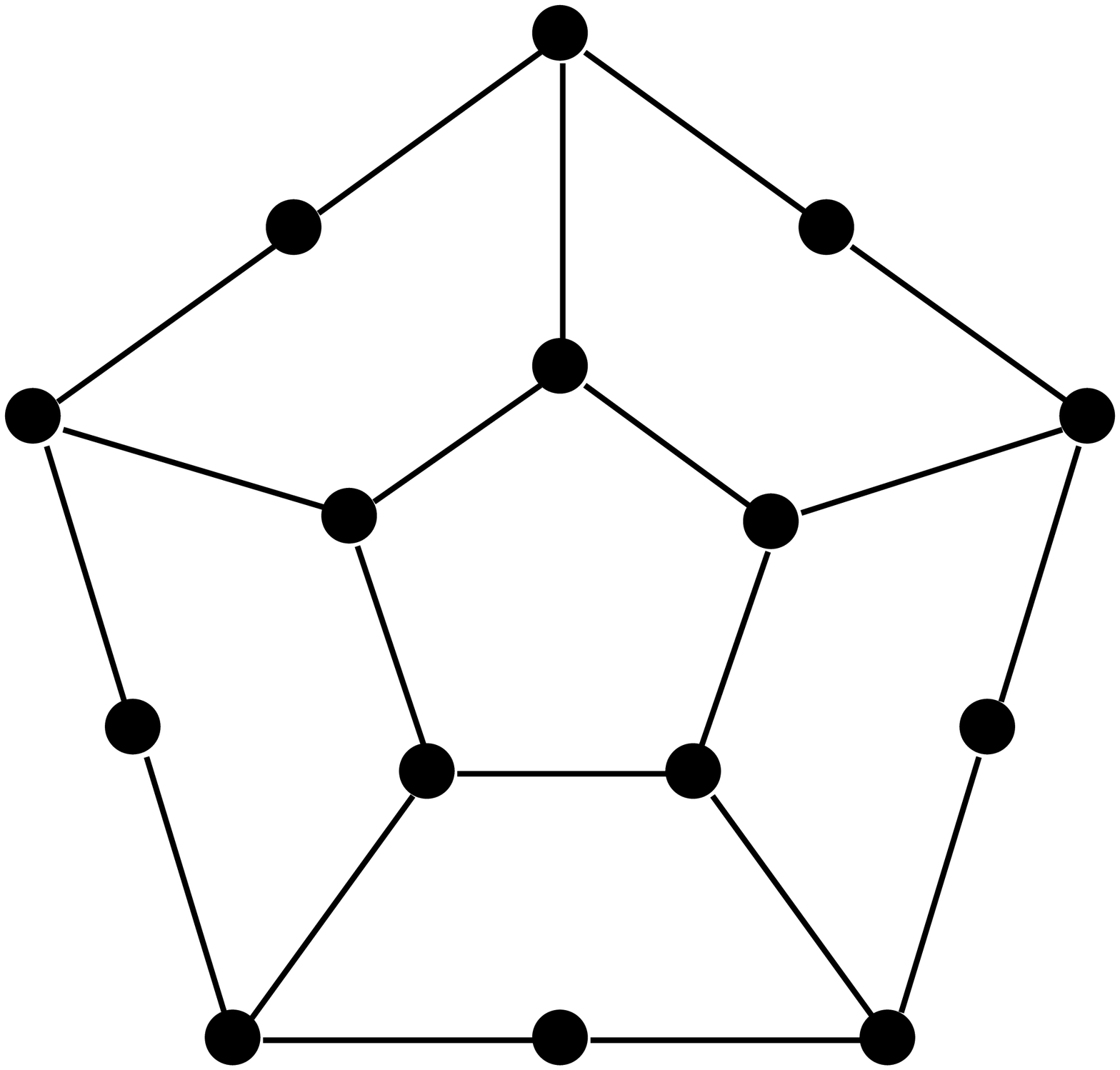, width=0.3\textwidth}
}}
\centerline{\hspace*{5mm}Dodecahedral graph \hspace*{1cm}Petersen graph
\hspace*{1cm}
Half-dodecaheral graph
}
\end{figure}

The paper is organized as follows. In the section 2, we will review
the definitions and some facts about Ricci curvature on graphs.
Theorem \ref{t1} will be proved in section 3. In the last section, we will discuss
general constructions using Cartesian
product, strong graph covering, and lattice graphs constructed by
 the root systems
of simple Lie algebras. Combining these methods, we constructed many
Ricci-flat graphs
with girth 3 or 4.

\section{Notations and Lemmas}
Let $G=(V,E)$ be a simple graph.
A probability distribution (over the vertex-set $V(G)$) is a mapping
$m\colon V \to [0,1]$ satisfying $\sum_{x\in V}m(x)=1$. Suppose that two
probability distributions $m_1$ and $m_2$ have finite support.
 A coupling between $m_1$ and $m_2$ is a mapping $A\colon V\times V\to [0,1]$
with finite support so that
$$\sum_{y\in V}A(x,y)=m_1(x) \mbox{ and } \sum_{x\in V}A(x,y)=m_2(y).$$

The transportation
distance between two probability distributions $m_1$ and $m_2$
is defined as follows.
\begin{equation}
  \label{eq:w1}
  W(m_1,m_2) = \inf_{A}\sum_{x,y\in V}A(x,y)d(x,y),
\end{equation}
where the infimum is taken over all couplings $A$ between $m_1$ and $m_2$.
A function $f$ over $G$ is $c$-Lipschitz if
$$|f(x)-f(y)|\leq c d(x,y)$$
for all $x,y\in V$.
By the duality theorem of a linear optimization problem,
the transportation distance can also be written as follows.
\begin{equation}
  \label{eq:w2}
  W(m_1,m_2) = \sup_{f}\sum_{x \in V}f(x)[m_1(x)-m_2(x)],
\end{equation}
where the supremum is taken over all $1$-Lipschitz function $f$.

For any vertex $x$, let $\Gamma(x)$ denote the
set of neighbors of $x$. I.e,
\begin{equation*}
  \Gamma(x)=\{v | vx \in E(G)\}.
\end{equation*}
Let $N(x)=\Gamma(x)\cup \{x\}$.

For any $\alpha\in [0,1]$ and any vertex $x$,
the probability measure $m_x^{\alpha}$ is defined as
\begin{equation}
  \label{eq:1}
  m^{\alpha}_x(v)=\left\{
    \begin{array}{ll}
      \alpha & \mbox{ if } v=x;\\
      \frac{1-\alpha}{d_x} & \mbox{ if } v\in \Gamma(x);\\
      0 & \mbox{ otherwise.}
    \end{array}\right.
\end{equation}

For any $x,y\in V$, we define $\alpha$-Ricci-curvature $\kappa_\alpha$ to be
\begin{equation}
  \label{eq:kalpha}
  \kappa_\alpha(x,y)=1-\frac{W(m^\alpha_x, m^\alpha_y)}{d(x,y)}.
\end{equation}

The Ricci curvature $\kappa(x,y)$ is defined as
$$\lim_{\alpha\to 1}\frac{\kappa_\alpha(x,y)}{1-\alpha}.$$

\begin{figure*}[hbt] 
 \begin{multicols}{2}
\centering
\unitlength=7.5mm
\begin{picture}(7,2)
\put(0,1){\line(1,0){6}}  
\put(1,0){\line(0,1){2}}  
\put(1,1){\circle*{0.2}}
\put(2,1){\circle*{0.2}}
\put(1,2){\circle*{0.2}}
\put(1,0){\circle*{0.2}}
\put(0,1){\circle*{0.2}}

\put(6,1){\circle*{0.2}}
\put(5,1){\circle*{0.2}}

\put(6,1){\line(1,2){0.447}}
\put(6,1){\line(1,-2){0.447}}
\put(6.447,1.894){\circle*{0.2}}
\put(6.447,0.106){\circle*{0.2}}
\put(0.75,0,65){$x$}
\put(5.75,0,65){$y$}

\put(1.05,1.15){$1-\epsilon$}
\put(1.1,1.75){$\frac{\epsilon}{d_x}$}
\put(1.1,0.25){$\frac{\epsilon}{d_x}$}
\put(0.1,0.65){$\frac{\epsilon}{d_x}$}
\put(2.1,0.65){$\frac{\epsilon}{d_x}$}

\put(6.2,0.65){$1-\epsilon$}
\put(5.2,0.65){$\frac{\epsilon}{d_y}$}
\put(6.547,1.75){$\frac{\epsilon}{d_y}$}
\put(6.547,0.25){$\frac{\epsilon}{d_y}$}
\end{picture} 
\caption{{\it The transportation distance between two $\epsilon$-balls
    in a graph.}
  \label{fig:gball}} \newpage
  \psfig{width=0.3\textwidth, file=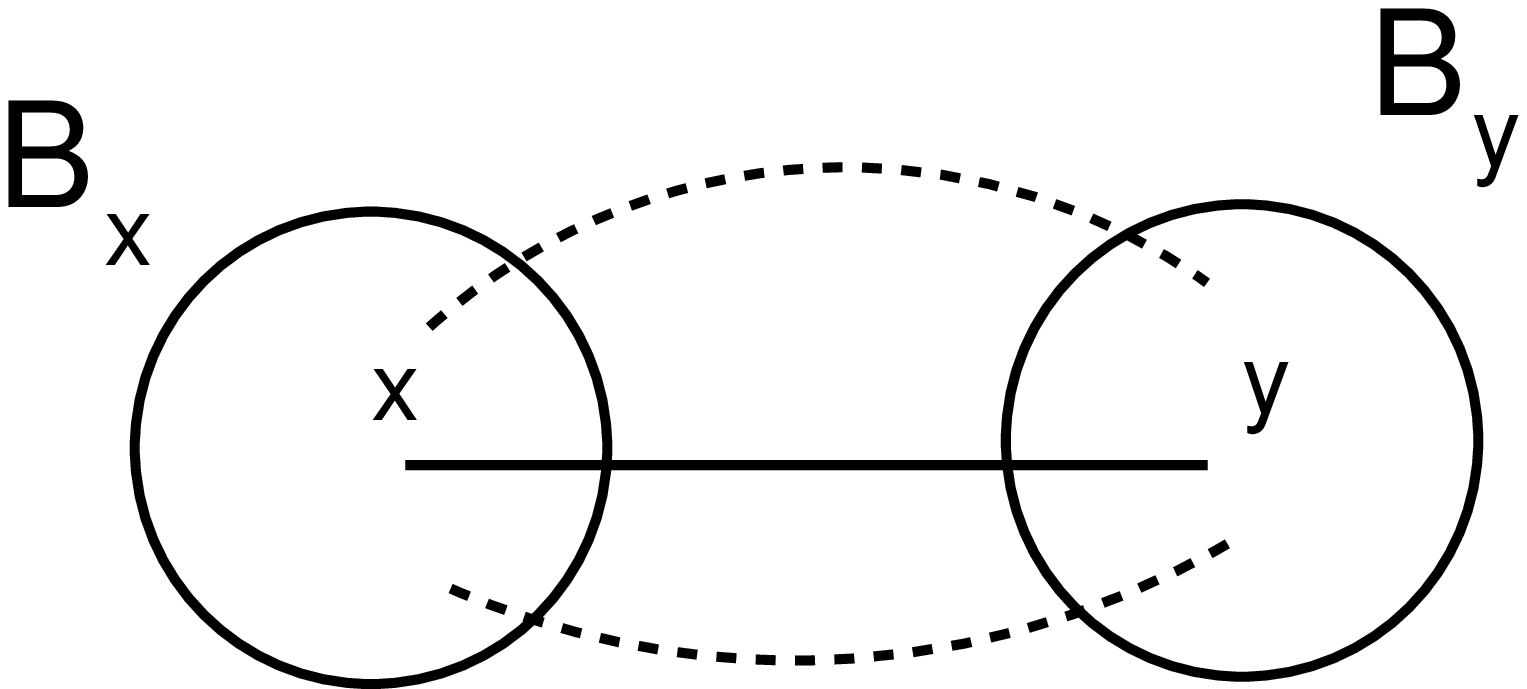}
\caption{{\it The average distance between two $\epsilon$-balls in
    a manifold.}
  \label{fig:mball}}
\end{multicols}
\end{figure*}

Let $\epsilon=1-\alpha$. Then $m^\alpha_x$ can be
viewed as the ball of radius $\epsilon$ and centered at $x$. Under
this setting, the Ricci curvature capture the approximation of the
transportation distance between two balls of radius $\epsilon$:
\begin{equation}
  \label{eq:mball}
W(m_x^{1-\epsilon}, m_y^{1-\epsilon})=(1-\epsilon \kappa(x,y)+o(\epsilon))d(x,y).  
\end{equation}
This is similar to the Ricci-curvature in the
differential manifolds $M$ of dimension $N$ (see Figures
\ref{fig:gball} and \ref{fig:mball}).
Let $x\in M$ and $v$ be a unit tangent vector at $x$. Let $y$ be a
point on the geodesic issuing from $v$, with $d(x,y)$ small enough. Then,
the average distance between two $\epsilon$-balls
centered at $x$ and $y$ in $M$ is  (see \cite{ollivier}):
\begin{equation}
  \label{eq:gball}
  d(x,y)\left(1-\frac{\epsilon^2Ric(v,v)}{2(N+2)}+O(\epsilon^3 +\epsilon^2d(x,y))\right).
\end{equation}

A graph $G$ is called {\it Ricci-flat} if $\kappa(x,y)=0$ for any edge $xy\in E(G)$.

\begin{lemma}\label{l1}
Suppose that an edge $xy$ in a graph $G$ is not in any $C_3$, $C_4$, or $C_5$.
Then $\kappa(x,y)=\frac{2}{d_x} + \frac{2}{d_y}-2$.
\end{lemma}

\begin{corollary}\label{c1}
Suppose that $x$ is a leaf-vertex (i.e., $d_x=1$). Let $y$ be the only neighbor of $x$.
Then $\kappa(x,y)>0$.
\end{corollary}

\begin{lemma}\label{l2}
Suppose that an edge $xy$ in a graph $G$ is not in any $C_3$ or $C_4$.
Then $\kappa(x,y)\leq \frac{1}{d_x} + \frac{2}{d_y}-1$.
  \end{lemma}
{\bf Proof:}
Define $f\colon N(x)\cup N(y)\to {\mathbb R}$ as follows.
$$f(u)=
\begin{cases}
  0 & \mbox{ if } u\in N(x) \setminus \{y\},\\
  1 & \mbox{ if } u=y,\\
  2 & \mbox{ if } u\in N(y)\setminus \{x\}.
\end{cases}
$$
Since $xy$ is not  in any $C_3$ or $C_4$, $f$ is a Lipschitz function.
We have
\begin{align*}
  W(m^\alpha_x, m^\alpha_y) &\geq \sum_{y\in V}f(u)[m_y^\alpha(u)-m_x^\alpha(u)]\\
& =\alpha +2(1-\alpha)\left(1-\frac{1}{d_y}\right)-\frac{1-\alpha}{d_x} \\
& = 2 - \alpha - (1-\alpha)\left(\frac{1}{d_x}+\frac{2}{d_y}\right).
\end{align*}
We have
$$\kappa(x,y)=\lim_{\alpha\to 1} \frac{1-W(m^\alpha_x, m^\alpha_y)}{1-\alpha}
\leq \frac{1}{d_x}+\frac{2}{d_y} -1.$$
The proof of Lemma is finished. \hfill $\square$

\begin{figure}[htbp]
\centerline{ {\psfig{figure=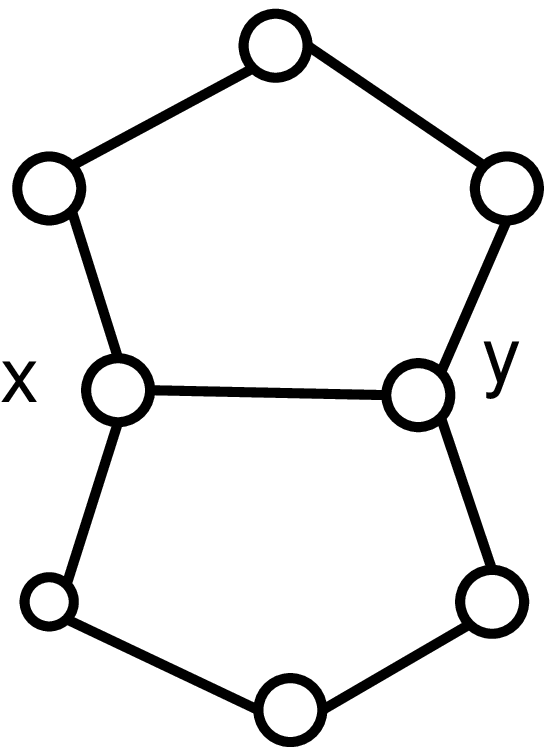, width=0.2\textwidth}}
\hfil \psfig{figure=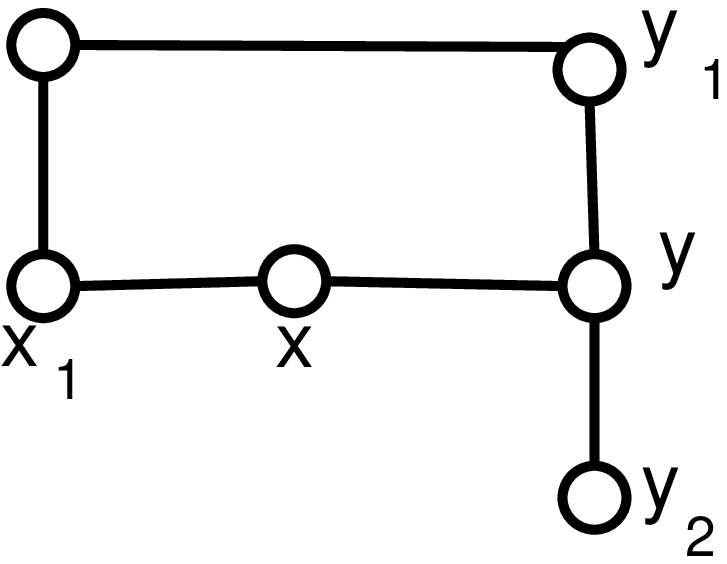, width=0.3\textwidth}
\hfil {\psfig{figure=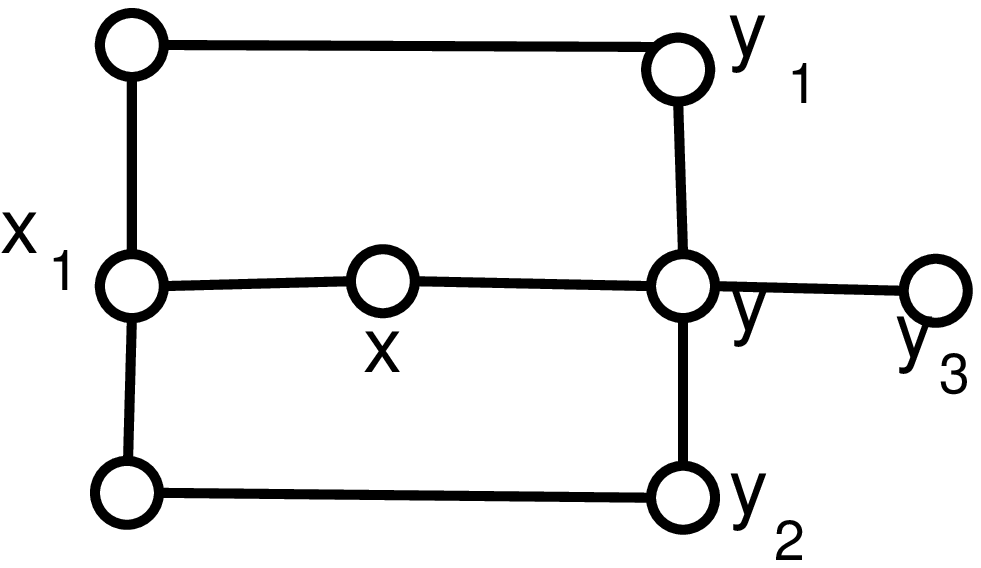, width=0.3\textwidth}
}}
\centerline{\hspace*{5mm}$d_x=d_y=3$ \hspace*{5mm} $d_x=2$ and $d_y=3$
\hspace*{1cm}
$d_x=2$ and $d_y=4$
}
\caption{Local structures in Lemma \ref{l3}.}\label{fig:3}
\end{figure}

\begin{lemma}\label{l3}
Suppose that an edge $xy$ in a graph $G$ is not in any $C_3$ or $C_4$.
Without loss of generality, we assume $d_x\leq d_y$.
If $\kappa(x,y)=0$, then one of the following statements holds (See Figure \ref{fig:3}).
\begin{enumerate}
\item $d_x=d_y=2$. In this case, $xy$ is not in any $C_5$.
\item $d_x=d_y=3$. In this case,  $xy$ is shared by two $C_5$s.
\item $d_x=2$  and $d_y=3$.  In this case,
let $x_1$ be the other neighbor of $x$ other than $y$.
Let $y_1,y_2$ be two neighbors of $y$ other than $x$.  Then
$\{d(x_1,y_1), d(x_1, y_2)\}=\{2,3\}$.
\item $d_x=2$ and $d_y=4$.  In this case, let $x_1$ be the other
  neighbor of $x$ other than $y$.  Let $y_1,y_2, y_3$ be three neighbors of
  $y$ other than $x$.  Then at least two of $y_1,y_2, y_3$ have distance 2 from $x$.
\end{enumerate}
\end{lemma}

{\bf Proof:} By Lemma \ref{l2}, we have
$$0\leq \frac{1}{d_x}+\frac{2}{d_y}-1.$$
Since $d_x\leq d_y$, we have $1\leq \frac{3}{d_x}$. By Corollary \ref{c1},
$d_x\geq = 2$. We have $d_x=2$, or $3$.
 If $d_x=2$, then $d_y=2,3,4$. If $d_x=3$, then $d_y=3$.

{\bf Case 1:} $d_x=d_y=2$ and $xy$ is in a $C_5$, then $\kappa(x,y)>0$. Contradiction!

{\bf Case 2:}  $d_x=d_y=3$ and $xy$ is not shared by  two $C_5$'s. Let $x_1,x_2$ be two other
neighbors of $x$ and $y_1, y_2$ are two other neighbors of $y$.
Without loss of generality, we can assume $d(x_1, y_1)\geq 3$ and $d(x_1,y_2)\geq 3$.
In this case, we can define a Lipschitz function $f\colon N(x)\cup N(y)\to {\mathbb R}$
as follows, $f(x_1)=-1$, $f(x_2)=f(x)=0$, $f(y)=1$, and $f(y_1)=f(y_2)=2$.
Using this function $f$, we can show $\kappa(x,y)<0$. Contradiction.

{\bf Case 3:} $d_x=2$ and $d_y=3$. Let $x_1$ be the other
neighbors of $x$ and $y_1, y_2$ are two other neighbors of $y$.
If $d(x_1,y_1)=d(x_1,y_2)=2$, then a similar calculation shows $\kappa(x,y)>0$.
 If $d(x_1,y_1)=d(x_1,y_2)=3$, then  $\kappa(x,y)=\frac{2}{2}+\frac{2}{3}-1<0$.
In both cases, we get contradiction.

{\bf Case 4:}  $d_x=2$ and $d_y=4$. Let $x_1$ be the other
neighbors of $x$ and $y_1, y_2, y_3$ are two other neighbors of $y$.
If $d(x,y_i)$ ($i=1,2,3$) have at least two $3$'s, then $\kappa(x,y)<0$.
Contradiction!

The proof of Lemma is finished. \hfill $\square$

\section{Proof of Theorem \ref{t1}}
{\bf Proof of Theorem \ref{t1}:}
Suppose $g(G)\geq 6$.  By Lemma \ref{l1}, for any edge $xy$, we have
$$0=\kappa(x,y)=\frac{2}{d_x} + \frac{2}{d_y}-2.$$
The only integer solution is $d_x=d_y=2$. Hence $G$ is $2$-regular graph.
Since $G$ is connected, it is either an infinite path or a cycle $C_n$.

Now we consider the case $g(G)=5$.  By Lemma \ref{l3}, every vertex in $G$
has degree $2$, $3$, or $4$.
We have the following three claims:

\noindent {\bf Claim a:}  $G$ has no vertex with degree $4$.

Suppose that $y$ is a vertex with degree $4$. By Lemma \ref{l3}, all its neighbors
must have degree $2$. Let $x$ be one of its neighbors and $x_1$ be the other
neighbor of $x$ other than $y$. The local structure is
shown in Figure \ref{fig:3}. We observe that $x_1$ have degree at least $3$.
If $d_{x_1}=3$, by Lemma \ref{l3}, $x_1x$ can not be shared by two $C_5$;
contradiction! If $d_{x_1}=4$, then all its neighbors have degree $2$.
Now we have two vertices of degree $2$ in a $C_5$ that are adjacent to each other.
This is a contradiction to Item 1 of Lemma \ref{l3}. Thus, we proved Claim a.

\noindent {\bf Claim b:} Any $C_5$ in $G$ can have at most one vertex with degree $2$.

\begin{figure}[hbt]
  \centering
\psfig{width=0.6\textwidth, file=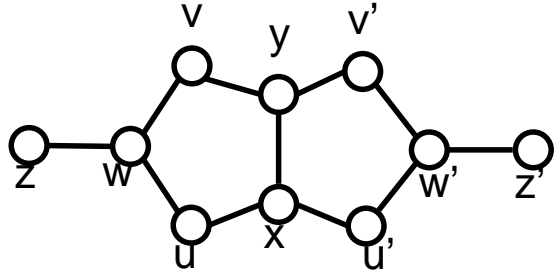}
  \caption{Two $C_5$s.}
  \label{fig:4}
\end{figure}

Suppose that $u$ and $v$ are two vertices of degree $2$
located in one $C_5$. By Item 1 of Lemma \ref{l3}, the vertices with degree $2$
are not adjacent to each other. In particular, the other three vertices on this $C_5$
have degree $3$. Let $x,y$ be two adjacent vertices of degree $3$  on a $C_5$.
Applying Lemma \ref{l3}, there is another $C_5$ passing through $xy$.
The other vertices are labeled as shown in Figure \ref{fig:4}.

If $d_{u'}=3$,  then $xu'$ is in two $C_5$. As result, $ww'$ must be connected.
Then $vw$ is in two $C_5$s. Contradiction to Item 3 of Lemma \ref{l3}.
Hence $d_{u'}=2$. Similarly, we have $d_{v'}=2$.  Note that  $d_w=d_{w'}=3$.
Note $zw$ is not in any $C_5$. Contradiction. Claim b is proved.

\noindent {\bf Claim c:} Suppose that $xy$ is an edge of $G$ with $d_x=2$ and $d_y=3$.
Let $y_1$ and $y_2$ be the other two neighbors of $y$ besides $x$. Then
one of $y_1$ and $y_2$ has degree two while the other one has degree three.

By Lemma \ref{l3}, every vertex of degree $2$ is in a $C_5$. This $C_5$ must contain $y$
and one of $y_1$ and $y_2$. Without loss of generality, assume that $y_1$ is in this $C_5$.
By Claim b, $y_1$ must have degree $3$.  By Lemma \ref{l3}, $yy_1$ is shared by
two $C_5$s. The second $C_5$ must pass through $y_2$. Since $xy$ is only in one $C_5$,
there is no $C_5$ containing $y_2$, $y$, and $x$. Thus, $d_{y_2}=2$.

If $G$ has no vertex with degree $2$, then $G$ is $3$-regular.
Every edge is in two $C_5$'s. $G$ can be embedded into a surface $S$.

If $S$ is an oriented surface with genus $s$, then Euler formula
gives
$$n-e+f=2-2s.$$
Here $e$ is the number of edges and $f$ is the number of faces.
Since $G$ is $3$-regular, we have $2e=3n$. Since every face is a $C_5$,
we have $2e=5f$. Thus
$$2-2s=n-e+f=\frac{2}{3}e- e \frac{2}{5}e=\frac{1}{15}e>0.$$
We must $s=0$. Hence $e=30$, $n=20$, and $f=10$. In this case,
$G$ is the dodecahedral graph.

If $S$ is a non-oriented surface with odd genus $s$, then Euler formula
gives
$$n-e+f=2-s.$$
Here $e$ is the number of edges and $f$ is the number of faces.
Since $G$ is $3$-regular, we have $2e=3n$. Since every face is a $C_5$,
we have $2e=5f$. Thus
$$2-s=n-e+f=\frac{2}{3}e- e \frac{2}{5}e=\frac{1}{15}e>0.$$
We must $s=1$. Hence $e=15$, $n=10$, and $f=6$. In this case,
$G$ is the Petersen graph.

\begin{figure}[hbt]
  \centering
\psfig{width=0.5\textwidth, file=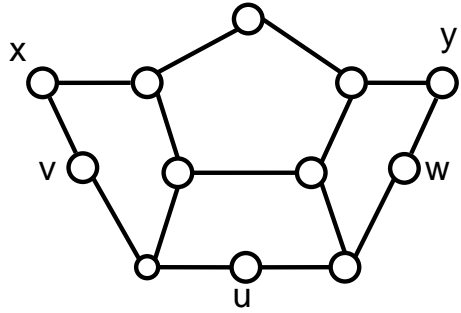}
  \caption{A partial graph of the half-dodecahedral graph.}
  \label{fig:5}
\end{figure}

If $G$ contains a vertex  $u$ with degree $2$, then $u$ is in a $C_5$.
The other four vertices have degree $3$. By Lemma 3, this forces
four $C_5$s as shown by Figure \ref{fig:5}.

By Claim c, we have $d_v=d_w=2$.
By Claim b, we have $d_x=d_y=3$. Now the graph continues to expand.
It finally results in the half-dodecahedral graph.
The proof of Theorem  is finished. \hfill $\square$

\section{More Ricci-flat graphs}
In this section, we will construct many Ricci-flat graphs
with girth $3$ or $4$. Given two graphs $G$ and $H$,
the Cartesian product (denoted by $G\square H$) is a graph
over the vertex set $V(G)\times V(H)$, where
two pairs $(u_1,v_1)$ and $(u_2,v_2)$ are connected if
``$u_1=u_2$ and $v_1v_2\in E(H)$'' or ``$u_1u_2\in E(G)$ and
$v_1=v_2$''. The following theorem is proved in \cite{LLY}.

{\bf Theorem (see \cite{LLY})}
{\it Suppose that $G$ is $d_G$-regular and $H$ is $d_H$-regular.
Then the Ricci curvature of $G\square H$ is given by
\begin{eqnarray}
  \label{eq:product1}
  \kappa^{G\square H}((u_1,v), (u_2,v))&=&\frac{d_G}{d_G+d_H} \kappa^G(u_1,u_2)\\
 \label{eq:product2}
  \kappa^{G\square H}((u,v_1), (u,v_2))&=&\frac{d_H}{d_G+d_H} \kappa^H(v_1,v_2).
\end{eqnarray}
Here $u\in V(G)$, $v\in V(H)$, $u_1u_2\in E(G)$, and $v_1v_2\in E(H)$.
}

\begin{corollary}
  If both $G$ and $H$ are Ricci-flat regular graphs, so is the
  Cartesian product graph $G\square H$.
\end{corollary}

Another contruction of Ricci-flat graphs is using a strong graph
covering. A graph $G$ is a {\em strong covering graph} of another graph $H$
if there is a surjective map $f\colon V(G)\to V(H)$ satisfying
for any edge $uv\in E(G)$, $f$ maps the induced graph of $G$ on
$\Gamma_G(u)\cup \Gamma_G(v)$
to the induced graph of $H$ on $\Gamma_H(f(u))\cup \Gamma_H(f(v))$
bijectively.
Note in the traditional
definition of covering graph, it only requires $f$ maps the induced
graph of $G$ on $\Gamma_G(u)$ to $H$ on $\Gamma_H(f(u))$ bijectively.
A strong covering graph $H$ is always a cover graph of $H$, not
vice versa. Observe that $\kappa(u,v)$ only depends on the induced
graph on $\Gamma(u)\cup \Gamma(v)$. We have the following corollary.

\begin{corollary}
  If $G$ is a strong covering graph of $H$, then
$G$ is Ricci-flat if and only if $H$ is Ricci-flat.
\end{corollary}

\begin{center}
  \includegraphics[width=4cm]{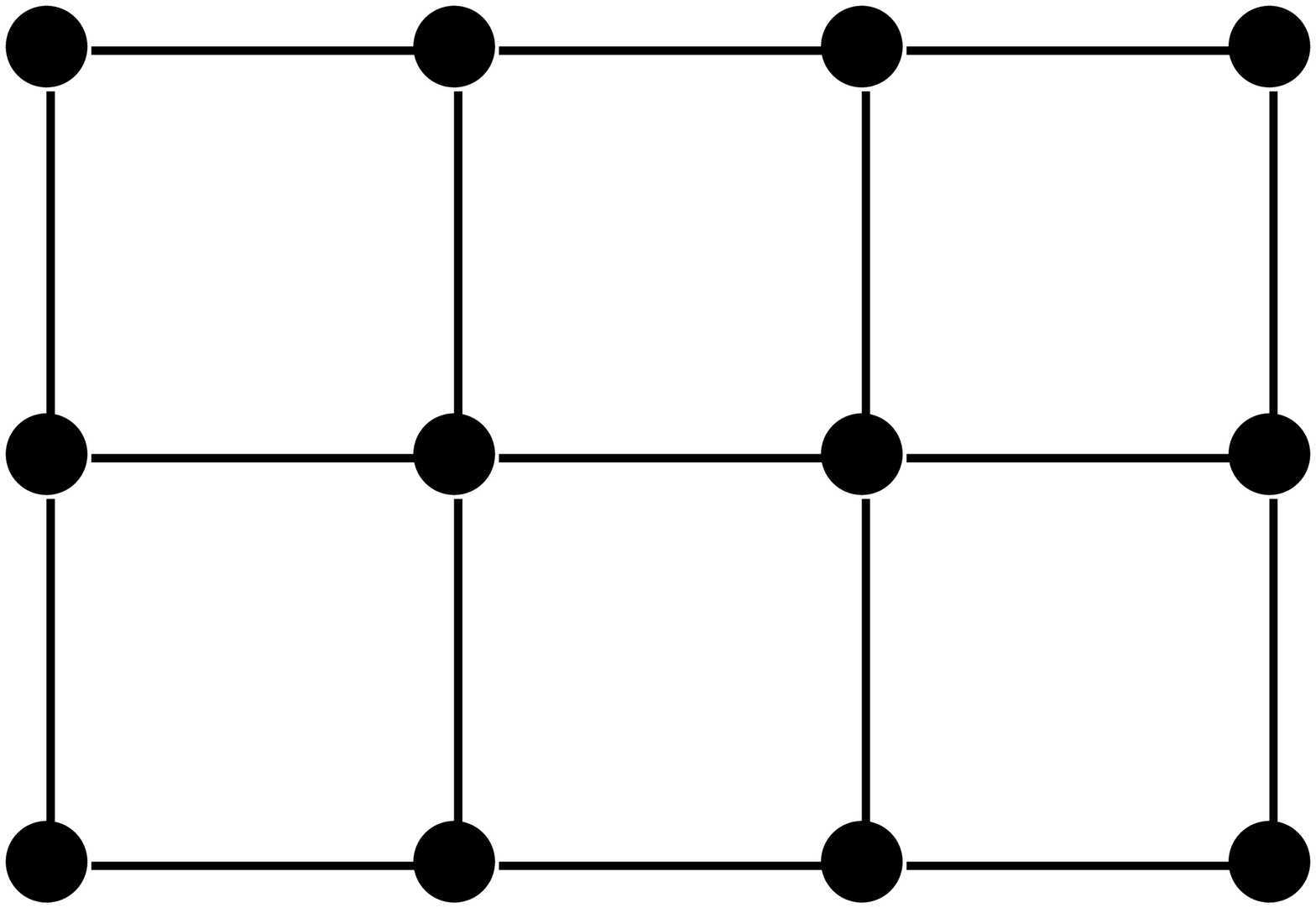}
\hfil
  \includegraphics[width=3cm]{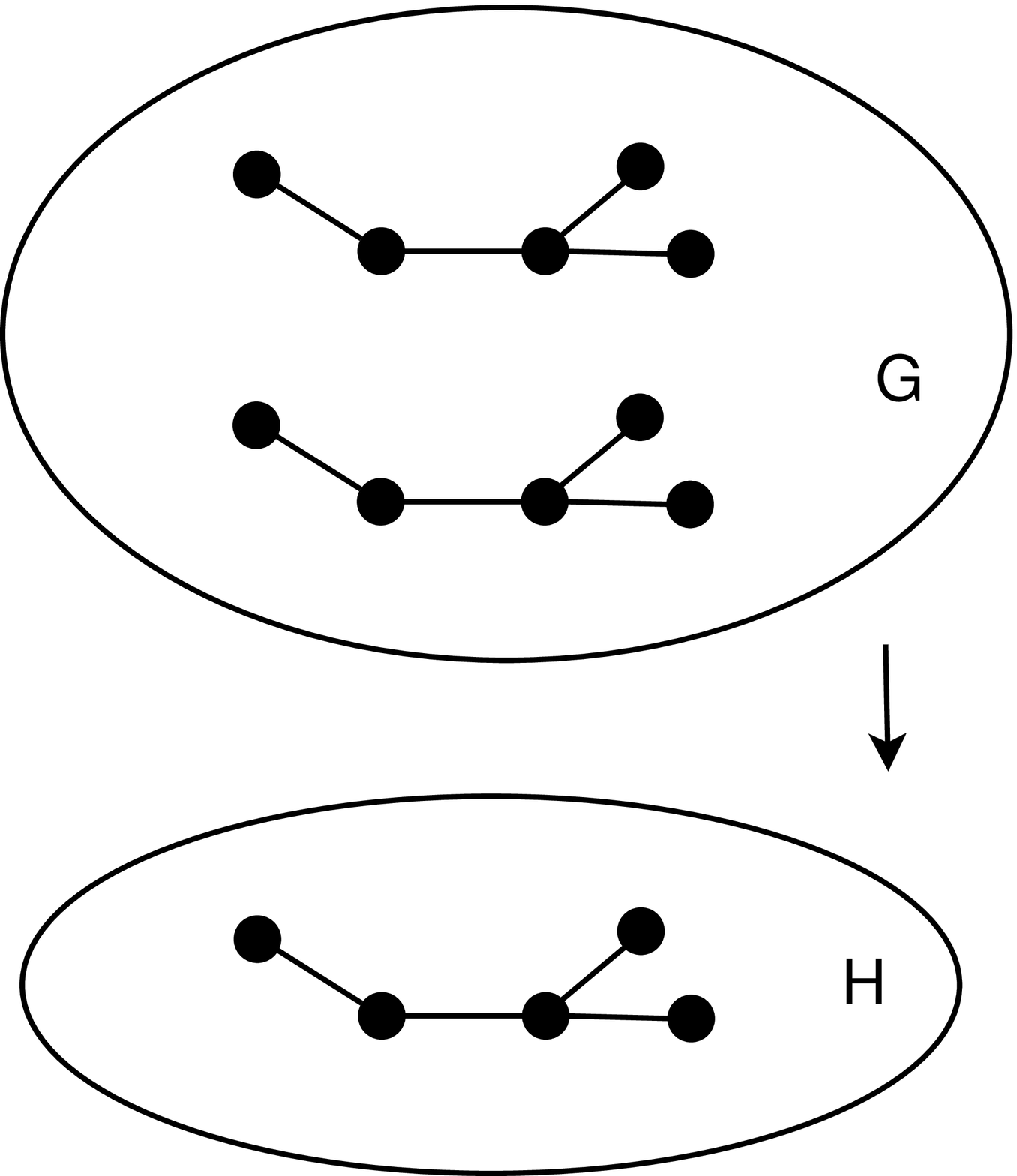}\\
\hspace*{1cm}
Cartesian Product \hspace*{2cm}
Strong covering graphs
\end{center}

A Ricci-flat graph $G$ is called {\em simple} if
$G$ is neither a Cartesian product of  Ricci-flat graphs
 nor a projection of some Ricci-flat graphs.


Many simple Ricci-flat graphs can be constructed using Cayley graphs.
Suppose that $G$ is a group and $S$ is a subset of $G$ with $-S=S$.
 The Cayley graph $\Cayley(G,S)$ is a graph with the vertex set of
 all elements in $G$ and where a pair $(x,y)$ forms an edge if
$x^{-1}y\in S$.

\begin{theorem} \label{cayleyinf}
Let $G$ be a torsion-free abelian group and $S$ be a finite subset of $G$
satisfying $S=-S$. Let $H:=\Cayley(G,S)$ be the Cayley graph.
Suppose that $H$ has following property
``for any $x\in V$, $s\in S$, and $n\geq 1$, the graph distance $d_H(x,x+ns)=n$.''
Then $H$ is Ricci-flat.
\end{theorem}
{\bf Proof:}
We first show $\kappa(x,y)\geq 0$ for any edge $xy$.
We use  addition ``+'' for the group operation of the abelian $G$.
Note that $H$ is vertex-transitive; thus it is regular with degree $d:=|S|$.
For $\alpha\in [0,1]$, we need construct a transportation $A$ that
moves the probability distribution $m^\alpha_x$ To $m^\alpha_y$ as
follows
$$
A(u,v)=
\begin{cases}
  \alpha & \mbox{ if } u=x \mbox{ and }v=y;\\
  \frac{1-\alpha}{d}  & \mbox{ if } u\in \Gamma_H(x) \mbox{ and }
  v-u=y-x;\\
  0 &\mbox{ otherwise.}
\end{cases}
$$
Observe that
$A$ simply moves every mass at vertex $u$ to $v:=u+ (y-x)$. Since $y-x\in
S$, $uv$ is always an edge. The cost of this transportation is
$\sum_{u,v\in G}A(u,v)d(u,v)=1$. By equations \eqref{eq:w1} and
\eqref{eq:kalpha}, we have $\kappa_\alpha(x,y)\geq 0$. Thus,
$\kappa(x,y)\geq 0$.

Now we show $\kappa(x,y)=0$. We will prove it by contradiction.
Suppose $\kappa(x,y)>0$ for some edge $xy\in E(H)$.
Let  $s:=y-x$ and $\epsilon:= \kappa(x,y)$. We have $s\in S$ and
$\epsilon>0$.
Since $\lim \limits_{\alpha\to
  1}\tfrac{\kappa_\alpha(x,y)}{1-\alpha}=\epsilon$,
there exists an $\alpha<1$ such that
$$\kappa_\alpha(x,y)>\frac{\epsilon}{2}(1-\alpha).$$

Choose an integer $n>\frac{4}{\epsilon}$.
Consider a path $x, x+s, x+2s, \ldots, x+ns$.
For any $g\in G$,  $g$ acts on $H$ as a graph homomorphism.
thus
$$\kappa(x,y)=\kappa(g+x,g+y).$$
In particular,  for $1\leq i\leq n$,
we have
$$\kappa_\alpha(x +(i-1)s,
x+is)=\kappa_\alpha(x,y)>\frac{\epsilon}{2}(1-\alpha).$$
From equation \eqref{eq:kalpha}, we get
$$W(m_{x+(i-1)s}^\alpha, m_{x+is}^\alpha)<
1-\frac{\epsilon}{2}(1-\alpha).$$
From the triangular inequality of the transportation
distance, we have
\begin{align*}
  d(x,x+ns)&= W(m^{1}_{x}, m^{1}_{x+ns}))\\
&\leq W(m^{1}_{x}, m^{\alpha}_{x}) +\sum_{i=1}^n W(m^{\alpha}_{x+(i-1)s},
m^{\alpha}_{x+is}) \\
&\hspace*{4mm} + W(m^{\alpha}_{x+ns}, m^{1}_{x+ns})\\
&< (1-\alpha) + n \left(1-\frac{\epsilon}{2}(1-\alpha)\right)
+(1-\alpha)\\
&=n+ (1-\alpha)\left(2- \frac{n\epsilon}{2}\right)\\
&<n.
\end{align*}
In the last step, we apply $2<\frac{n\epsilon}{2}$ by our choice of $n$.
Contradiction!
\hfill $\square$

Now we  apply Theorem \ref{cayleyinf} to some special lattice graphs
whose vertices are the root lattices in some Euclidean spaces.

\begin{theorem} \label{rootlattice}
Let $R$ be a root system of type $A_n$, $B_n$, $C_n$, $D_n$, $F_4$,
$E_6$, $E_7$, or $E_8$ (except $G_2$). Let $G$ be the root lattice,
which is viewed as the abelian additive
group generated by the roots in $R$.
 Let $S\subset R$ satisfying $-S=S$. Then the Cayley graph
 $\Cayley(G,S)$  is  Ricci-flat.
\end{theorem}

{\bf Proof:} Write $H:=\Cayley(G,S)$. It suffices to show the
conditions of Theorem \ref{cayleyinf} are satisfied. We only need
to verify the following condition:
``{\it For any $x\in G$ and $s\in S$, $d_H(x,x+ns)=n$.}''

By the triangular inequality of the graph distance, we have 
$$d_H(x,x+ns)\leq \sum_{i=1}^n d_H(x+(i-1)s, x+is)\leq n.$$
We only to verify $d_H(x,x+ns)\geq n$.
We will compare the graph distance to the Euclidean distance as $G$
is embedded as a root lattice. There are two cases.

{\bf Case 1:}
 The root system $A_n$, $D_n$, $E_6$, $E_7$,
and $E_8$ have only one type of length, say $l$. For any $x\in G$ and $s\in S$,
 For any $x,y\in G$,  $d_H(x,y)$ be the graph distance while
$d(x,y)$ be the Euclidean distance. Then we have
$$d(x,y)\leq d_H(x,y)l.$$
For any $s\in S$ and integer $n$, we have
$$d_H(x,x+ns)\geq \frac{1}{l}d(x,x+ns)= n \frac{d(0,s)}{l}=n.$$
 By Theorem \ref{cayleyinf}, $\Cayley(G,S)$ is
Ricci-flat.

{\bf Case 2:}
Now we consider the root system of type $B_n$, $C_n$, and $F_4$.
Those root systems have two types of lengths. The ratio of the length of
the longer root to one of the shorter root is $\sqrt{2}$.
 In addition, the angle $\theta$ formed between a longer root $u$ and a shorter
root $v$ is either $\frac{\pi}{4}$, $\frac{\pi}{2}$, or
$\frac{3\pi}{4}$.
We have
\begin{equation}
  \label{eq:uv}
\langle u, v\rangle=\|u\| \|v\|\cos \theta= \epsilon \|v\|^2,
\end{equation}
where $\epsilon\in \{-1,0,1\}$.
Let $l$ denote the length of the shorter root. Then the length of
the longer root is $\sqrt{2}l$. Equation \eqref{eq:uv} implies
for any shorter root $v$ and any root $u$ (whether short or long):
$$\langle u, v\rangle\leq l^2.$$

If $s$ is a  longer root, then we have
$$d_H(x,x+ns)\geq \frac{1}{\sqrt{2}l}d(x,x+ns)= n \frac{d(0,s)}{\sqrt{2}l}=n.$$
We are done.

If $s$ is a shorter root and $d_H(x, x+ns)=k<n$, the we can find
$s_1,s_2,\ldots, s_k \in S$ such that
$ns=s_1+s_2+\cdots+s_k$.
Taking an inner product with $s$, we get
$$nl^2=\langle ns, s\rangle = \sum_{i=1}^k <s_i,s> \leq kl^2.$$
We get $k\geq n$. Contradiction!

By Theorem \ref{cayleyinf}, $\Cayley(G,S)$ is always
Ricci-flat. \hfill $\square$.

The rank of a root system is the dimension of the underlined Euclidean
space. The only rank-1 root system is $A_1$, which consists of
two roots $\alpha$ and $-\alpha$. The Cayley graph generated by this
root system is the infinite path.

The rank-2 root systems $R$ are $A_2$, $B_2=C_2$, $D_2=A_1+A_1$, and
$G_2$. Since $G_2$ is excluded from Theorem \ref{rootlattice},
Here we get three lattice graph $\Cayley(G, R)$. (See Figures
\ref{fig:A1A1}, \ref{fig:A2}, \ref{fig:B2}, where the root systems
are drawn in red.)
We didn't get more
lattice graphs by taking a subset
$S\subseteq R$. For example, for the root system $B_2$ and $S$
consisting of $3$ pairs of roots, we get a ``skewed'' drawing of
the graph generated by $A_2$. However, for the root system $R$ of rank more
than 2, we do get more Ricciflat lattice graphs by taking proper subset $S\subsetneq R$.

\begin{figure}[hbt]

\begin{center}
\tikzstyle{every node}=[circle, draw, fill=black!50, inner sep=0pt,
minimum width=4pt]
\begin{tikzpicture}[thick,scale=0.4]
\pgfmathsetmacro{\n}{7};
\draw \foreach \x in {0,1,...,\n} {
      \foreach \y in {0,1,...,\n} {
          (\x,\y) node {} -- (\x+1,\y)
         (\x,\y)  -- (\x,\y+1)
      }
  };

 \draw \foreach \x in {0,1,...,\n} {
    (\x,\n+1) node {} -- (\x+1,\n+1)
   (\n+1,\x) node {} -- (\n+1,\x+1)
 };
\draw (\n+1,\n+1) node {};
\draw[very thick, red, ->, fill=red] (4,4) -- (4,5);
\draw[very thick, red, ->] (4,4) -- (5,4);
\draw[very thick, red, ->] (4,4) -- (4,3);
\draw[very thick, red, ->] (4,4) -- (3,4);
\draw[very thick, red, ->] (4,4) node[fill=red] {};
  \end{tikzpicture}\quad
\begin{tikzpicture}[thick,scale=0.4]
\pgfmathsetmacro{\n}{8};
\foreach \j in {0,1,...,\n} {
   \pgfmathsetmacro{\h}{\n-\j};
   \foreach \i in {0,...,\h} {
      \pgfmathsetmacro{\x}{\i+0.5*\j};
      \pgfmathsetmacro{\y}{0.866*\j};
        \draw  (\x,\y) node {} -- (\x+1,\y)
          (\x,\y)  --  (\x+0.5,\y+0.866)
          (\x+1,\y) -- (\x+0.5,\y+0.866);
       }
  }
\pgfmathsetmacro{\m}{\n+1};
\foreach \i in {0,1,...,\m} {
      \pgfmathsetmacro{\x}{\m-\i*0.5};
      \pgfmathsetmacro{\y}{0.866*\i};
     \draw (\x,\y) node {};
}
      \pgfmathsetmacro{\x}{4.5};
      \pgfmathsetmacro{\y}{0.866*3};
\draw[very thick, red, ->]  (\x,\y)  --  (\x+1,\y);
\draw[very thick, red, ->]  (\x,\y)  -- (\x+0.5,\y+0.866);
\draw[very thick, red, ->]  (\x,\y)  --  (\x-0.5,\y+0.866);
\draw[very thick, red, ->] (\x,\y)  --  (\x-1,\y);
\draw[very thick, red, ->] (\x,\y)  --  (\x-0.5,\y-0.866);
\draw[very thick, red, ->] (\x,\y)  --  (\x+0.5,\y-0.866);
\draw (\x,\y) node[fill=red] {};

  \end{tikzpicture}\quad
\begin{tikzpicture}[thick,scale=0.4]
\pgfmathsetmacro{\n}{7};
\draw \foreach \x in {0,1,...,\n} {
    \foreach \y in {0,1,...,\n} {
          (\x,\y) node {} -- (\x+1,\y)
          (\x,\y)  -- (\x,\y+1)
          (\x,\y)  -- (\x+1,\y+1)
          (\x+1,\y)  -- (\x,\y+1)
      }
  };
 \draw \foreach \x in {0,1,...,\n} {
    (\x,\n+1) node {} -- (\x+1,\n+1)
   (\n+1,\x) node {} -- (\n+1,\x+1)
 };
\draw (\n+1,\n+1) node {};
\draw[very thick, red, ->] (4,4) -- (4,5);
\draw[very thick, red, ->] (4,4) -- (5,4);
\draw[very thick, red, ->] (4,4) -- (4,3);
\draw[very thick, red, ->] (4,4) -- (3,4);
\draw[very thick, red, ->] (4,4) -- (5,5);
\draw[very thick, red, ->] (4,4) -- (5,3);
\draw[very thick, red, ->] (4,4) -- (3,5);
\draw[very thick, red, ->] (4,4) -- (3,3);
\draw[very thick, red, ->] (4,4) node[fill=red] {};
  \end{tikzpicture}
\end{center}
  \begin{multicols}{3}
\caption{{\it Lattice graph with root system $A_1+A_1$.}
  \label{fig:A1A1}} \newpage
\caption{{\it Lattice graph with root system $A_2$.}
  \label{fig:A2}} \newpage
\caption{{\it Lattice graph with root system $B_2$.}
  \label{fig:B2}}
\end{multicols}
\end{figure}

Combining with the strong covering graph, we get the following corollary.
\begin{corollary}
  Let $H:=\Cayley(G,S)$ be the lattice graph constructed in Theorem
  \ref{rootlattice}. Let $G'$ be a subgroup of $G$ satisfying
``any two elements in $G'$ has a graph-distance at least $6$ in $H$.''
Define the quotient graph $H/G'$ with vertex set $G/G'$ and edge set
$u+G'\sim v+G'$ if $u-v\in SG'$. Then $H/G'$ is Ricci-flat.
\end{corollary}

We get many examples of finite Ricci-flat graphs, which can be viewed
as the discrete analogue of torus.

\end{document}